\documentstyle[amssymb,11pt]{article} 
\newfont{\gothic}{eufm10 scaled 1100} 
\def \goths{{\hbox{\gothic S}}} 
\newfont{\smallgothic}{eufm10 scaled 700} 
\def \smallgoths{{\hbox{\smallgothic S}}} 
\newcommand{\gap}{{\Gamma_{1,p}}}

\newcommand{\alptilde}{{\tilde{\calA}_{1,p}}}
\newcommand{\alpptilde}{{\tilde{\calA}_{{p^2}}}}
\newcommand{\alpast}{{\calA_{1,p}^\ast}}
\newcommand{\openalp}{{\open{\calA_{1,p}}}}
\newcommand{\alppast}{{\calA_{{p^2}}^\ast}}
\newcommand{\openalpp}{{\open{\calA_{{p^2}}}}}
\newcommand{\alpr}{{\calA_{1,p}''}}
\newcommand{\alpbar}{{\overline{\calA}_{1,p}}}
\newcommand{\Do}{D_{{\ell_0}}} 
\newcommand{\Dab}{D_{{\ell_{a,b}}}}
\newcommand{\Symp}{{\, \mbox{Sp}}} 
\newcommand{\calA}{{\cal A}}
\newcommand{\calL}{{\cal L}}
\newcommand{\calN}{{\cal N}} 
\newcommand{\calS}{{\goths}}
\newcommand{\calSk}{{{{\cal S}}_k}} 
\newcommand{\oka}{{\cal O}}
\newcommand{\CC}{{\Bbb C}} 
\newcommand{\PP}{{\Bbb P}}
\newcommand{\QQ}{{\Bbb Q}} 
\newcommand{\ZZ}{{\Bbb Z}}
\newcommand{\Eins}{\mbox{$1 \hspace{-4pt}1$}} 
\newcommand{\Thm}[1]{{\bf #1. $\:$}} 
\newcommand{\proof}{{\em Proof. $\:$}}
\newcommand{\open}[1]{#1^\circ} 
\newcommand{\ebew}{\hfill$\Box$ \par}
\newcommand{\trace}{\mbox{\rm trace}\,}
\newcommand{\Text}[1]{\qquad \mbox{#1} \quad} 
\newtheorem{thm}{Theorem}[section] 
\newtheorem{prop}[thm]{Proposition}
\newtheorem{lemma}[thm]{Lemma} 
\newtheorem{defi}[thm]{Definition}
\newtheorem{cor}[thm]{Corollary}
\newtheorem{rk}[thm]{Remark}

\begin{document} \title{Invariants of Moduli Spaces of Abelian
Surfaces} \author{J\"org Zintl} \date{02.02.99} \maketitle

Our objects of study are compactified moduli spaces of {$(1,p)$-polarized}\linebreak 
abelian surfaces with level structure of canonical type. They are obtained
from quotients of Siegel space $\calS_2$ by a certain paramodular group
$\gap$. The construction of the compactification has been described by
Hulek, Kahn and Weintraub in \cite{HKW}. Our aim is to determine invariants of a desingularized model $\alptilde$. This could
serve as a test case for similarly constructed moduli spaces. In
particular we are interested in numerical invariants of $\alptilde$.  To
avoid combinatorial complications we will restrict our attention to $p \ge
5$ and prime.

The aim of this article is twofold. First we want to give a short account
of the results obtained in \cite{Z}. These include a geometric description
of the canonical divisor of $\alptilde$ as well as formulae for its self
intersection number and the Euler number of the moduli space. We will give
brief indications of proofs. Because computations in the toroidal compactification as in \cite{HKW} soon become technical and lenghty, the full arguments cannot be repeated here. Instead we try to emphasize the main ideas of the approach. 
It was partially stimulated by a paper by Yamazaki \cite{Y}, who considered the case of moduli spaces of abelian
surfaces with level-$n$-structure $\calA_n^\ast$. The main theorem of \cite{Z} is

\Thm{Theorem}{\em The canonical divisor of $\alptilde$ equals \[
K_\alptilde = 3 L - D - \frac{1}{2} R - \frac{1}{2} E . \] }

Here $L$ denotes the $\QQ$-divisor associated to modular forms, $D$ the
boundary divisor and $R$ and $E$ arise from the ramification of the
quotient map $\calS_2 \rightarrow \gap \setminus \calS_2$.

In the second part we will compute a formula for the second Chern number
 $c_1(\alptilde). c_2(\alptilde)$. To do this we count modular forms on $\alptilde$ applying a fixed point theorem of Atiyah and Singer \cite{AS} and compare this to a calculation using the theorem of Riemann--Roch.
This makes use of intersection numbers on $\calA_n^\ast$ determined by Yamazaki \cite{Y}. Thus we obtain as combined result: 

\Thm{Theorem}{\em Put $c_i := c_i(\alptilde)$. The Chern numbers of the
moduli space $\alptilde$ are \[ \begin{array}{lcl} {c_1}^3 &=& -
\frac{p^2-1}{960} (9p^3 - 360 p^2 + 1519 p + 3000) \\[2mm]
 c_1 . c_2 &=& -
\frac{p^2-1}{240} (p-13) (p^2 - 17 p + 90) \\[2mm] c_3 &=& - \frac{p^2-1}{1440} (
p^3 + 431 p - 8760) .  \end{array} \] }

Since $c_1.c_2$ is a birational invariant for three-dimensional varieties
this implies that there are no unirational moduli spaces other than the
known ones for $p = 5, 7$ and $11$. See also \cite{HM,MS,GP} and \cite{Gr}. 

We conclude with the deduction  of a
polynomial expression for the dimension of spaces of modular forms and cusp forms with respect to $\gap$. 

I would like to thank Professor K. Hulek very much for his advice and
encouragement during the completion of \cite{Z}.

\section*{Part I: R\'esum\'e}
\setcounter{section}{1}
For us a moduli space will be a quotient
$\Gamma \setminus \calS_2$, where $\calS_2$ denotes the Siegel space of
degree 2, that is \[ \calS_2 := \{ \tau \in \: \mbox{Sym}_2(\CC) : \:
\mbox{Im} (\tau) \: \: \mbox{positive definite} \: \} \] and $\Gamma
\subset \Symp_4(\QQ)$ an arithmetic subgroup.  As usual the action is by fractional linear transformations. 
For the rest of this paper
$p$ will denote an integer, which for technical reasons will be assumed to
be prime and $p \ge 5$. Put \[ \gap := \left\{ \gamma \in \Symp_4(\ZZ) :
\gamma - \Eins \in \left( \begin{array}{rrrr} \ZZ & \ZZ& \ZZ & p\ZZ\\
p\ZZ&p\ZZ&p\ZZ&p^2\ZZ\\ \ZZ&\ZZ&\ZZ&p\ZZ\\ \ZZ&\ZZ&\ZZ&p\ZZ \end{array}
\right) \right\} .\]

\begin{defi}\em The quotient $\openalp := \gap \setminus \calS_2$ is called
the {\em (open) moduli space of $(1,p)$--polarized abelian surfaces with level
structure of canonical type.}
\end{defi}

This is a quasi-projective variety of dimension 3 and it has only finite
quotient singularities. Note that the quotient map $\pi : \calS_2
\rightarrow \openalp$ is ramified over two disjoint surfaces in
$\openalp$, which are called {\em Humbert surfaces} $\open{H_1}$ and
$\open{H_2}$. For more details see [HKW 1,2]. 

Let $\Gamma$ denote some arithmetic subgroup of $ \Symp_4(\QQ)$. We denote
by $M_k(\Gamma)$ the space of modular forms of weight $k$ with respect to
$\Gamma$. By definition $f \in M_k(\Gamma)$ is a holomorphic function $f :
\calS_2 \rightarrow \CC$ satisfying \[ f(\gamma \tau) = \: \mbox{det} (C
\tau + D)^k f(\tau) \] for all $\gamma = {A B \choose C D}  \in \Gamma$ and all $\tau \in \calS_2$. For $k$
sufficiently divisible $ \mbox{det} (C \tau + D)^k$ is a factor of
automorphy defining a line bundle $\calL_{{\Gamma \setminus \smallgoths_2}}^k$
on $\Gamma \setminus \calS_2$. Modular forms of sufficiently large weight
$k$ define an embedding of $\Gamma \setminus \calS_2$ into projective
space 
\[ \iota_k : \Gamma \setminus \calS_2 \; \: \hookrightarrow \: \: \PP^N = \PP
H^0(\Gamma \setminus \calS_2, \calL_{{\Gamma \setminus \smallgoths_2}}^k) \] with $\iota_k^\ast \, \oka_{{\PP^N}}(1)
= \calL_{{\Gamma \setminus \smallgoths_2}}^k$. The projective closure
$\overline{\Gamma \setminus \calS_2}$ of $\Gamma \setminus \calS_2$ in
$\PP^N$ is isomorphic to the Satake compactification  \cite{Sa} of $\Gamma \setminus
\calS_2$. In particular the line bundle $\calL_{{\Gamma \setminus
\smallgoths_2}}^k$ can be extended to a line bundle on $\overline{\Gamma
\setminus \calS_2}$. Observe that the boundary $\delta := \overline{\Gamma
\setminus \calS_2} - \Gamma \setminus \calS_2$ is of codimension two, and in general 
$\overline{\Gamma \setminus \calS_2}$ is highly singular along $\delta$.

One can use the theory of toroidal compactifications as described in
\cite{AMRT} to obtain less singular compactifications of $\Gamma \setminus
\calS_2$. Such a construction was carried out in great detail by Hulek,
Kahn and Weintraub in \cite{HKW} in the case of $\Gamma = \gap$. Because
of combinatorial complications they restrict their attention to the case
$p$ prime. They define the Igusa compactification \[ \alpast := \openalp
\cup D \] by attaching a boundary divisor $D$ consisting of one {\em central
boundary component} $\Do$ whose normalization is isomorphic to the Kummer
modular surface $K(p)$ and $\frac{p^2-1}{2}$ {\em peripheral boundary
components} $\Dab$ isomorphic to $K(1)$. For the indexing of the boundary components we refer to \cite{HKW}.
 There exists a natural morphism
onto the Satake compactification \[ \alpast \: \: \rightarrow \: \: \alpbar :=
\overline{\gap \setminus \calS_2} .\]

\Thm{Remark} There are a number of results on the compactified moduli
space $\alpast$. For example Hulek and Sankaran showed in \cite{HS2} that $\alpast$ is
simply connected. It is a variety of general type for $p \ge 37$ by
\cite{HS, GH} and rational if $p = 5, 7$ \cite{HM,MS}. 

In \cite{HKWsing} a full description of the singularities of $\alpast$
was given. 

\begin{prop} If $p \ge 5$ the singular locus of $\alpast$ consists of two
smooth curves $C_1$ and $C_2$ isomorphic to $X(p)$, both contained in
$H_1$, which is the closure of $\open{H}_1$, as well as two additional isolated singular points $Q_{a,b}^{(3)}$
and $Q_{a,b}^{(4)}$ lying in each peripheral boundary component $\Dab$.
For $i=1,2$ the curve $C_i$ intersects the boundary in precisely one point
$Q_{a,b}^{(i)}$ of each peripheral boundary component $\Dab$. \end{prop}

We choose a smooth desingularisation $\alptilde$ obtained by blowing up
 $\alpast$ along $C_1$ and $C_2$ and resolving the isolated singularities
in a straightforward way. See \cite{Z} for details. Denote by
$\calL_{{\alpast}}$ and $\calL_\alptilde$ the lifting of the $\QQ$-line
bundle $\calL_{{\alpbar}}$ to $\alpast$ and $\alptilde$ respectively.
$\calL_\alptilde^{k}$ is a line bundle if $k$ is divisible by $12$.

\begin{defi}\em  Let $m = (m',m'') \in \ZZ^2 \times \ZZ^2$. We define the {\em
theta constant of characteristic $m$} on $\calS_2$ by \[ \Theta_m(\tau) :=
\sum_{{q \in \ZZ^2}} \: \mbox{\rm exp}\, \left( 2 \pi i \left[ \frac{1}{2} (q
+ \frac{m'}{2}) \tau (q + \frac{m'}{2})^t + (q+
\frac{m'}{2})(\frac{m''}{2})^t \right] \right) \] for $\tau \in \calS_2$.
There are precisely $10$ different theta constants which do not vanish
identically on $\calS_2$. We put \[ \Theta^2(\tau) := \prod_{{\Theta_m
\not\equiv 0}} \Theta_m^2 (\tau) .\] \end{defi}

This holomorphic function is ubiquitous in mathematics. It is known to be a
modular form of weight 10 with respect to $\Symp_4(\ZZ)$. Up to scalar
multiplication it is even the unique cusp form of weight $10$. See for example 
\cite{I}. 

For this reason $\Theta^{12}$ defines a section $s_\Theta$ in
$\calL_{{\openalp}}^{60}$ which can be lifted and extended to a section in
$\calL_\alptilde^{60}$. By a result of Hammond \cite{Ham} the modular form
$\Theta^2$ vanishes precisely on \[ \Delta := \bigcup_{{\gamma \in \;
\mbox{\scriptsize Sp}_4(\ZZ)}} \gamma \Delta_0 \quad \:  \mbox{where} \:\:
\Delta_0 := \left\{ \left( \begin{array}{cc} \tau_1&\tau_2\\ \tau_2&\tau_3
\end{array} \right) \in \calS_2: \tau_2 = 0 \right\}. \] The image of
$\Delta$ in $\openalp$ decomposes into two disjoint irreducible
components. One of them is the Humbert surface $\open{H_1}$, the other one
a smooth surface $\open{T_1}$. Since $\gap$ is not a normal subgroup of
$\Symp_4(\ZZ)$, the component $\open{T_1}$ is {\em not} a translate of
$\open{H_1}$ under a group action.

Denote the closed liftings of $\open{H_1}$, $\open{H_2}$ and $\open{T_1}$
in 
 $\alptilde$ by $H_1$, $H_2$ and $T_1$. Let $E^{(1)}$ and $E^{(2)}$ denote
the exceptional divisors arising from the blow up along $C_1$ and $C_2$.
Over each of the isolated singular points $Q_{a,b}^{(3)}$ on
a peripheral boundary component $\Dab$ there lies an exceptional divisor $E_{a,b}^{(3)}$, and we write $E_{a,b}^{(4)}$ and $E_{a,b}^{(5)}$ for the two irreducible
components of the exceptional divisors of the isolated singular points $Q_{a,b}^{(4)}$ as in \cite{Z}. 
We denote the strict transforms of the boundary divisors in $\alptilde$ again by the same symbols. 
Using this we
can give a geometrical description of a divisor $60 L$ associated to
$\calL_\alptilde^{60}$, i.e. associated to modular forms of weight $60$. 

\begin{thm} The associated divisor $60 L$ on $\alptilde$ equals \[
\begin{array}{lcl} 60 L &=& 6 H_1 + 12 T_1 + 6 \Do + 6
\sum\limits_{{\ell_{a,b}}} c_{a,b} \Dab\\ && + 3 E^{(1)} + 2 E^{(2)} +
\sum\limits_{{\ell_{a,b}}} 3 c_{a,b} E^{(3)}_{a,b} +
\sum\limits_{{\ell_{a,b}}} c_{a,b} ( 2 E^{(4)}_{a,b} + 4 E^{(5)}_{a,b} )
\end{array} \] where $ c_{a,b} = \left\{ \begin{array}{lcl} p^2 & ,& a
\equiv 0 \: \mbox{\rm mod} \; p\\ 1& ,& \mbox{otherwise.} \end{array}
\right. $
\end{thm}
Peripheral boundary components with $a \equiv 0$ mod $p$ are called {\em standard} components, and all others are called {\em nonstandard} components.  

For the proof see \cite[Satz 4.2.9]{Z}. There we study the zero divisor of $s_\Theta$
on $ \openalp$. Using the explicit description of the toroidal
compactification of $\openalp$ in \cite{HKW} one can write down $s_\Theta$
in local coordinates in an open neighbourhood of the boundary components. This shows immediately that the section extends 
over the boundary. Finally one reads off the coefficients of the zero
divisor on $\alptilde$ as stated above.

Along the same lines one can give a description of the canonical divisor. 

\begin{thm}\label{1}
 The canonical divisor $K = K_\alptilde$ of $\alptilde$ equals
\[ K = 3 L - D - \frac{1}{2} R - \frac{1}{2} E \] where \[
\begin{array}{lcll} D &=& \Do + \sum_{{\ell_{a,b}}} \Dab &\mbox{(boundary
divisor)}\\ R&=& H_1 + H_2&\mbox{(ramification divisor)}\\ E&=&
\frac{1}{2} E^{(1)} + E^{(2)} + \frac{1}{2} \sum_{{\ell_{a,b}}}
E_{a,b}^{(3)} &\mbox{(exceptional contribution)}. \end{array} \] \end{thm}

Again, this was shown in \cite[Thm. 4.3.4]{Z}. One has to consider coordinates $\tau = {\tau_1 \: 
\tau_2 \choose \tau_2 \: \tau_3}$ on $\calS_2$. Then \[ \omega(\tau) :=
\Theta^{12}(\tau) \cdot (d \tau_1 \wedge d\tau_2 \wedge d\tau_3)^{20} \]
is a $\Symp_4(\ZZ)$-invariant form on $\calS_2$. In particular it defines
a section of \[ K_\alptilde^{20} | \openalp - (\open{H}_1 \cup \open{H}_2) . \] As
before we study the extension of the section to all of $\alptilde$. This
time we get a meromorphic section with zeroes and poles as given in the
theorem.

\Thm{Remark} Compare this to Yamazaki's \cite{Y} result. He considered the
moduli space $\calA_n^\ast$ of abelian surfaces with level-$n$-structures.
This is a toroidal compactification of the quotient $\calA_n^\circ :=
\Gamma_2(n) \setminus \calS_2$, where $\Gamma_2(n)$ is the principal
congruence subgroup of level $n$.
 Since this group is torsion free there is no ramification (and hence no
singular locus), and he finds \[ K_{{\calA_n^\ast}} = 3 L - D \] using the
appropriate definition of $L$ and $D$.

In Theorem \ref{1} we have an expression for $K_\alptilde$ in geometrical
terms: It involves the Humbert surfaces $H_1$ and $H_2$, the boundary divisors $\Do$
and $\Dab$, exceptional divisors and the $\QQ$-divisor $L$ associated to
modular forms. So there is a straightforward way to compute the self
intersection number of the canonical divisor. We use the identity \[
A.B.C = A|C . B|C \] for equivalence classes of smooth hypersurfaces $A, B$ and $C$. This reduces the problem to intersection theory on surfaces. In \cite{Z}
we gave a complete account of the groups of numerical equivalence classes
of divisors of all hypersurfaces involved. The computation of the relevant intersection tables is not difficult, but it requires some amount of preparations and notations. Hence we only quote the results of the calculations in tables 1,2 and 3 using the notations of theorem \ref{1} and defining $\kappa:= p^2-1$. 

{\em Table 1:} {Intersection products of type $A.B.C$ with $C = R$}
\[\begin{array}{r|rrrr}
\cdot& L&R&D&E\\[1pt]
\hline
L& \frac{7p}{144} \kappa &-\frac{19p}{144} \kappa& \frac{p+1}{6} \kappa&\frac{p}{16}\kappa\\[1pt]
R& & \frac{4p-51}{12}\kappa&- \frac{p+9}{3}\kappa & -\frac{p-2}{8}\kappa\\[1pt]
D&&&\frac{5}{4} \kappa &\kappa\\[1pt]
E&&&&- \frac{1}{4} \kappa
\end{array} \]

{\em Table 2:} {Intersection products of type $A.B.C$ with $C = E$}
\[ \begin{array}{r|rrrr}
\cdot&L&R&D&E\\
\hline
L&0&\frac{p}{16}\kappa&0&-\frac{7p}{48}\kappa\\[1pt]
R&&-\frac{p-2}{8}\kappa & \kappa & - \frac{1}{4} \kappa\\[1pt]
D&&&\frac{1}{4} \kappa&-2 \kappa\\[1pt]
E&&&&\frac{19p+6}{24}\kappa
\end{array} \]

{\em Table 3:} {Intersection products of type $A.B.D$ including $L$ and $D$ only}
\[ \begin{array}{r|rr}
\cdot& L&D\\
\hline
L&0& - \frac{p^2+1}{24} \kappa\\[1pt]
D&& - \frac{11p+18}{12} \kappa
\end{array} \]

We also know the self intersection number of $L$, which is $L^3 = \frac{p(p^4-1)}{2880}$.
Adding up all 
intersection numbers with appropriate (rational) coefficients one obtains
the following result. 

\begin{prop} \label{2} {\rm \cite[Thm. 5.1.34]{Z}}
The self intersection number of $K_\alptilde$ equals \[
{K_\alptilde}^3 = \frac{p^2-1}{960} (9p^3 - 360 p^2 + 1519 p + 3000) .\]
\end{prop}

As a by-product of the study of the hypersurface $\Do$ one finds a lower
bound for the Picard number of the moduli space $\alptilde$. 

\begin{prop} {\rm \cite[Satz 6.1.1]{Z}} 
The Picard number of $\alptilde$ is bounded from below by the following inequality:
\[ \rho(\alptilde) \ge 2 p^2 + 4 .\] 
\end{prop}

To see this one has to look at all divisors on $\alptilde$ we know so far. The tables of intersection numbers given in \cite{Z} show how many of them can possibly be independent. Most of them, precisely $2p^2+4$, are. 

Finally we quote a formula for the Euler number from \cite[Thm. 5.2.6]{Z}. 

\begin{prop}\label{3}
The Euler number of $\alptilde$ equals \[ \mbox{\rm e}
(\alptilde) = - \frac{p^2-1}{1440} ( p^3 + 431 p - 8760) .\] \end{prop}

To verify this note that $\Gamma_2(p^2)$ is a torsion free normal subgroup of $\gap$. This induces a
covering of open moduli spaces \[ q : \openalpp \rightarrow
 \openalp .\] The ramification behavior of $q$ has been well
described, see for example \cite{HKWsing}. Since the Euler number of
$\openalpp  = \Gamma_2(p^2) \setminus \calS_2$ is given by a
formula of Harder \cite{Harder} we easily obtain e$(\openalp)$. By
\cite{HKW} the boundary components and their mutual intersection behavior
are known. Adding up all contributions gives the proposition. 

Observe that the action of the quotient group $\Lambda := \gap / \Gamma_2(p^2)$ extends to Igusa's compactification $\alppast$ of $\openalpp$. In
\cite[section 6.3]{Z} we showed that there is a morphism \[ \rho : \alptilde \rightarrow
\Lambda \setminus \alppast \] as well as a $\Lambda$-invariant blow up $\alpptilde$ of
$\alppast$ to make the following diagramm commutative.  \[
\begin{array}{ccccc} &\alppast &\leftarrow &\alpptilde\\
\phi&\downarrow&&\downarrow &\tilde{\Phi}\\ &\Lambda \setminus \alppast&
\leftarrow & \alptilde \end{array} \] In fact $\rho$ is an isomorphism
outside the corank 2 boundary components of $\alptilde$. The quotient map $\phi$ is of degree $d := [\gap : \Gamma_2(p^2)] = p^{13}(p^2-1)$. The morphism
$\tilde{\Phi}$ is generically finite of the same degree, but it is not
induced by the group action of $\Lambda$.

\section*{Part II: The second Chern number} 
\setcounter{section}{2}
\setcounter{thm}{0}

In \cite{Z} we made a
conjecture on a formula for the Chern number $c_1 . c_2$ where $c_i := c_i(\alptilde)$ for $i=1,2$. The proof we
shall give here is a modification of the approach proposed there. By theorem \ref{1}  we get
\[ c_1 . c_2 = - ( 3 L - D - \frac{1}{2} R - \frac{1}{2} E) . c_2 \]
Apart from $c_2.L$ almost all numbers needed were computed in \cite {Z}, see table 4 further below.

The idea of the approach is the following. Let $S_k(\Gamma)$ denote the complex vector space of cusp forms of weight $k$ with respect to some arithmetic subgroup $\Gamma \subset \Symp_4(\QQ)$. We want to compute the dimension of $S_k(\gap)$ in two different ways.  To do this we will always assume $k \equiv 0$ mod $12$ so that $\oka_{\alptilde}(kL-D)$ is a line bundle.  Because of
\[ S_k(\gap) = H^0(\alptilde, \oka_{\alptilde}( k L - D) ) \]
we may apply the theorem of Riemann--Roch and vanishing theorems. This computation involves $c_2 . L$. On the other hand we can compare this to a dimension formula we obtain by exploiting the relation between $\alptilde$ and Yamazaki's space $\alppast$. 

Yamazaki computed in \cite{Y} a formula for the dimension of $S_k(\Gamma_2(p^2) )$ as a polynomial in $k$ and $p$.
\[ \begin{array}{rcl}
\dim S_k(\Gamma_2(p^2))\hspace{5mm} &=&\\[2mm]
= \frac{(p^2-1)^2}{2^8 3^3 5} \cdot [ \: \: \:  & k^3  \cdot &(2 p^{16} + 2 p^{14})\\
+&k^2  \cdot & (-9 p^{16} - 9 p^{14})\\
+&k^1 \cdot & (13 p^{16} + 13 p^{14} -120 p^{12} - 120 p^{10} )\\
+&k^0 \cdot & ( -6 p^{16} - 6 p^{14} + 180 p^{12} + 540  p^{10} + 360  p^8 ) \: \: ]
\end{array} \]

Note that there is a misprint in Yamazaki's paper \cite{Y}, page 39, formula $(v)$ as compared to his theorem 5.

By abuse of notation we denote the divisors corresponding to $L$ and $D$ on $\alppast$ by the same letters. Now the dimension formula above is  an application of the theorem of Riemann--Roch and the following vanishing theorem.

\begin{prop}\label{6}
Let $k \ge 4$ and $k \equiv 0$ mod $12$. Then
\[ \begin{array}{llcl}
&H^q(\alppast, \oka_{\alppast}(kL-D) ) &=& 0\\
\mbox{and}&H^q(\alptilde, \oka_{\alptilde}(kL-D) ) &=& 0
\end{array} \]
for all $q > 0$.
\end{prop}

\proof
For the first statement see \cite{Y}. 
Because of $k \equiv 0$ mod $12$ the $\QQ$-divisor $kL$ on $\alptilde$ is actually integral, and we can define an integral divisor
\[ M_k := (kL - D) - K_{\alptilde} = (k-3) L + \frac{1}{2} R + \frac{1}{2} E \]
and a $\QQ$-divisor
\[ M_k' := M_k - \frac{1}{2} R - \frac{1}{2} E . \]
Note that $M_k' = (k-3) L$, hence $M_k'$ is nef and big as a $\QQ$-divisor, with integral part $[M_k'] = M_k$. The support of the non-integral part $M_k' - [M_k'] = -\frac{1}{2} R - \frac{1}{2} E$ is a divisor with normal crossings. From the vanishing theorem of Kawamata and Viehweg \cite{K} we get
\[ H^q(\alptilde, \oka_{\alptilde}([-M_k']) ) = 0 \]
for $q < 3$ and therefore
\[ H^q(\alptilde, \oka_{\alptilde}(kL-D)) = H^q(\alptilde, \oka_{\alptilde}(K_{\alptilde} + M_k)) = 0 \]
for $q > 0$ by Serre duality.
\ebew

In the same manner as for $\oka_{\alppast}(kL-D)$ we can apply the theorem of Riemann--Roch to the line bundle $\oka_\alptilde(k L-D)$. Using the intersection numbers from tables 1,2 and 3  we get
\[ \dim S_k(\gap) = \hspace*{10cm} \]
\[ \begin{array}{crrll}
 &=&  &\frac{1}{12} &    (kL-D)(kL-D-K_\alptilde)(2kL-2D-K_\alptilde)\\
&&+& \frac{1}{12} &(kL-D).c_2 \: + \: \frac{1}{24} \: c_1.c_2 \\[2mm] 
&=&  \frac{p^2-1}{2^8 3^3 5} \cdot [ \: \: \:  & k^3  \cdot & (2 p^{3} + 2 p^{})\\
&&+&k^2 \cdot &  (-9 p^{3} + 201 p^{})\\
&&+&k^1 \cdot & ( 9  p^{3} - 120 p^{2} - 1481  p^{} -1080  ) \\
&&+&k^0 \cdot & ( \ldots ) \: ]\\[2mm]
&&+ \, \frac{1}{12} & k^1 \cdot &  L . c_2 \\[2mm]
&&+ \, \frac{1}{24} & k^0 \cdot &  c_1.c_2\, .
\end{array} \]

In particular $\dim S_k(\gap) = \dim H^0(\alptilde,\oka_\alptilde(kL-D))$ is a polynomial of degree $3$ in $k$. But there is a different approach to compute this polynomial. Remember that $\Gamma_2(p^2)$ is a normal subgroup of $\gap$, with quotient group $
\Lambda$ of order $d = p^{13}(p^2-1)$. 
Following a method of Hirzebruch \cite{Hi} we identify $S_k(\gap)$ with the $\Lambda$--invariant modular forms of $S_k(\Gamma_2(p^2))$ and we use the identity
\[ \begin{array}{lcl} 
d \cdot \dim S_k(\gap) &=& d \cdot \dim S_k^\Lambda(\Gamma_{2}(p^2))\\
&=& \sum\limits_{\gamma \in \Lambda} \trace( \gamma^\ast | S_k(\Gamma_2(p^2) )) .
\end{array} \]
This leads to the equation
\[ (\ast) \quad d \cdot \dim S_k(\gap) - \dim S_k(\Gamma_2(p^2)) = \sum_{\stackrel{\gamma \in \Lambda}{\gamma \neq \Eins}} \trace(\gamma^\ast | S_k(\Gamma_2(p^2))) .\]
We know the left hand side of this equation from our considerations above. Looking at the right hand side it will be possible to compute a formula for $c_2.L$. Note that is is enough to determine the coefficient of $k^1$ in this polynomial identity. 
Using this method the leading coefficient of the polynomial $\dim S_k(\gap)$ had been determined in \cite{HS}.

Put $\calSk := \oka_\alppast(k L - D )$. Since $H^q(\alppast, \calSk) = 0$ for $k \ge 4$, $k \equiv 0$ mod $12$ and all $q > 0$ by proposition \ref{6}, the fixed point theorem 
of Atiyah and Singer \cite{AS} shows
\[ \trace(\gamma^\ast | S_k(\Gamma_2(p^2))) = \sum_i (-1)^i \: \trace( \gamma^\ast | H^i(\alppast, \calSk) )  \]
\[ = \left\{ \frac{\mbox{ch}\,(\calSk | X^\gamma)(\gamma) \,\cdot\,  \prod\limits_\Theta U^\Theta(N^\gamma(\Theta)) \,\cdot\, \mbox{Td}\,(X^\gamma)}{\det(1-\gamma | (N^\gamma)^\ast) } \right\} [X^\gamma] .\]
Here $X^\gamma$ denotes the fixed set of $\gamma \in \Lambda$ in $\alppast$, Td$(X^\gamma)$ the Todd class of $X^\gamma$, $N^\gamma$ its normal bundle and $N^\gamma(\Theta)$ the subbundle of $N^\gamma$ such that $\gamma$ acts on $N^\gamma(\Theta)$ as $e^{i \Theta}$. 
Furthermore
\[ U^\Theta(N^\gamma(\Theta)) = \prod_\alpha \left( \frac{1- e^{{-x_\alpha - i \Theta}}}{1 - e^{i \Theta}} \right)^{-1} \Text{if} c(N^\gamma(\Theta)) = \prod_\alpha (1 + x_\alpha) .\]
Since $\calSk = \oka_\alppast(kL-D)$ is a line bundle we have
\[ (\dagger) \qquad \quad \mbox{ch}\,(\calSk | X^\gamma)(\gamma) = \chi(\gamma) \, \mbox{ch}\,(\calSk | X^\gamma)  \]
i.e. the usual Chern character times a character $\chi$ of the stabilizing group of $X^\gamma$. Note that the group $\Lambda$ acts trivially on $\oka_\alppast(kL)$, so $\chi(\gamma)$ is determined by the action of $\gamma$ on $\oka_\alppast(-D)| X^\gamma$. In all but one of our cases $\chi(\gamma)$ will be equal to $1$.  

The action of $\Lambda$ on $\alppast$ has been described in \cite{Z}, see also \cite{HKWsing}. Since we are only interested in the coefficients of nontrivial powers of $k$, we don't have to worry about fixed sets $X^\gamma$, on which $L | X^\gamma$ is trivial. In particular we may ignore isolated fixed points on the boundary and fixed curves on the corank 2 boundary components. 

Fixed sets of codimension 1 are the preimages $H_i^\ast$ of the Humbert surfaces $H_i$ fixed by involutions $\nu_i$, for $i = 1,2$, as well as the preimages $\Do^\ast$ and $\Dab^\ast$ of the central boundary component $\Do$ and the nonstandard peripheral 
boundary components $\Dab$. In each of the latter cases the stabilizer is a cyclic group of order $p^2$.

In codimension 2 contributing fixed sets are the preimages $C_i^\ast$ of the curves $C_i$ contained in $H_1$, with $i=1,2$. If $i=1$ then the stabilizer is generated by an element $\zeta_1$ of order $4$, if $i=2$ then by an element $\zeta_2$ of order $6$.
 We must also consider the curves $B_i^\ast$  of intersection of the preimages of $H_i$ for $i=1,2$ with those boundary components, which contribute to the ramification. In each case the stablizer is generated by the stabilizers of the intersecting ramification divisors and has order $2p^2$. 

We will determine each contribution separately.

1) Let us consider the case codim$\, X^\gamma = 1$ first. Since $N^\gamma(\Theta) = N^\gamma$ one computes
\[ \begin{array}{l}
 U^\Theta(N^\gamma(\Theta)) \quad = \\[2mm]
\quad = \left( \frac{1 - e^{i \Theta}}{1 - e^{-i \Theta}} \right) \left( 1 - \frac{e^{-i \Theta}}{1 - e^{-i \Theta}} \,  c_1(N^\gamma) + \frac{1}{2} \frac{e^{-i \Theta} (1+ e^{-i\Theta})}{(1 - e^{-i\Theta})^2} \, c_1(N^\gamma)^2 \right) 
\end{array} \]
as well as
\[ \begin{array}{lcl}
 \mbox{\rm ch}\, (\calSk | X^\gamma) &=& 1 + c_1(\calSk | X^\gamma) + \frac{1}{2} c_1(\calSk | X^\gamma)^2 \\
&=& 1 + (k L-D) |X^\gamma + \frac{1}{2} ((kL-D)|X^\gamma)^2 
\end{array} \]
and
\[ \begin{array}{lcl}
\mbox{\rm Td}\, (X^\gamma) &=& 1 - \frac{1}{2} K_{{X^\gamma}} + \frac{1}{12} ( K_{{X^\gamma}}^2 + c_2(X^\gamma) ) . 
\end{array} \]
Because of $\det (1 - \gamma | (N^\gamma)^\ast) = 1 - e^{i \Theta}$ and using $K_{{X^\gamma}} = (3 L - D + X^\gamma) | X^\gamma$ by the adjunction formula, we finally find
\[ \trace(\gamma^\ast | S_k(\Gamma_2(p^2))) = \hspace*{8cm}\]
\[ \begin{array}{rrl}
=  \frac{1}{1-e^{-i \Theta}} [\: & k^2 \cdot&  \frac{1}{2} (L|X^\gamma)^2 \\
+ & k^1 \cdot& ( - \frac{3}{2} (L | X^\gamma)^2 - \frac{1}{2} L D X^\gamma - \frac{1}{2} L (X^\gamma)^2 - \frac{e^{-i\Theta}}{1-e^{-i\Theta}} L (X^\gamma)^2 )\\
+ &  k^0 \cdot& ( \ldots ) \: ] .
\end{array} \]

\underline{Case 1.a:} Let $\gamma = \nu_1$ and $X^\gamma = H_1^\ast$ with $e^{i\Theta} = -1$. In this case we find
\[ \begin{array}{lrl}
t_1 &:=& \trace(\nu_1^\ast|S_k(\Gamma_2(p^2))) \\[1mm]
&=& \frac{1}{2} \, [ k^2 \cdot \frac{1}{2} L^2
 H_1^\ast + k^1 \cdot (- \frac{3}{2} L^2 H_1^\ast - \frac{1}{2} L D H_1^\ast) + \ldots \: ].
\end{array} \]
Let $H_1^{\ast,c}$ be an irreducible component of $H_1^\ast$. They are mutually disjoint, and each   $H_1^{\ast,c}$ is isomorphic to $X(p^2) \times X(p^2)$.
By Yamazaki's results
\[ \begin{array}{lcl}
L^2 H_1^{\ast,c} &=& \frac{p^8(p^2-1)^2}{12 \cdot 24} \\[10pt]
L (H_1^{\ast,c})^2 &=& - \frac{p^8(p^2-1)^2}{12 \cdot 24} .
\end{array} \]
Since the natural map 
\[ H_1^{\ast,c} \rightarrow H_1 \cong X(1) \times X(p) \]
has degree $\frac{p^7(p^2-1)}{2}$ and the quotient map $\phi$ is ramified of order $2$ above $H_1$, we conclude that there are 
\[ \frac{p^{13}(p^2-1)}{2} : \left( \frac{p^7(p^2-1)}{2} \right)^{-1} = p^6 \]
components of $H_1^\ast$, each of them meeting $p^2(p^2-1)$ boundary components isomorphic to $S(p^2)$ along a section isomorphic to $X(p^2)$. Note that 
\[ \mbox{deg}\, L | X(p^2) = \frac{p^4(p^2-1)}{24} .\] 
Therefore we compute 
\[ L D H_1^{\ast,c} = p^2(p^2-1) \cdot \frac{p^4(p^2-1)}{24} . \]
Hence we get 
\[ \begin{array}{lcl}
t_1 &=&  k^2 \cdot \frac{p^{14}(p^2-1)^2}{2^7 3^2} + k^1 \cdot ( - \frac{3}{4} \frac{p^{14}(p^2-1)^2}{2^5 3^2}  - \frac{1}{4}   \frac{p^{12}(p^2-1)^2}{2^3 3}) + \ldots \\[8pt]
&=& \frac{(p^2-1)^2}{2^8 3^3 5} ( k^2 \cdot 30 p^{14} +  k^1 \cdot (- 90 p^{14} - 360 p^{12}) + \ldots \: ) .
\end{array} \] 

\underline{Case 1.b:} Let  $\gamma = \nu_2$ and $X^\gamma = H_2^\ast$ with $e^{i \Theta} = -1$.
Note that for an irreducible component $H_2^{\ast,c}$ of $H_2^\ast$ there is a $6:1$ cover
\[ H_2^{\ast,c} \rightarrow H_1^{\ast,c} \]
ramified $2 : 1$ along the intersection of $H_2^{\ast,c}$ with the boundary of $\alppast$.  
Using this and the identity $H_2^{\ast,c} | H_2^{\ast,c} = K_{{H_2^{\ast,c}}} - (3L-D)|H_2^{\ast,c}$ one computes
\[ \begin{array}{lcl}
L^2 H_2^{\ast} &=& 6 \cdot L^2 H_1^\ast \\
L (H_2^\ast)^2 &=& 6 \cdot  L (H_1^\ast)^2
\end{array} \]
and
\[ L D H_2^{\ast,c} = \quad 3 \cdot L D H_1^{\ast,c} .\]
Hence the trace formula gives
\[ \begin{array}{lrl}
t_2 &:=& \trace(\nu_2^\ast | S_k(\Gamma_2(p^2))) \\
&=& \frac{(p^2-1)^2}{2^8 3^3 5} ( k^2 \cdot 180 p^{14} + k^1 \cdot (-540 p^{14} - 1080 p^{12}) + \ldots \: ) .
\end{array}\]

\underline{Case 1.c:} The boundary components $D^{\ast,c}$ with $e^{i \Theta} =e^{{2\pi i/p^2}} =: \xi$,
where $D^{\ast,c}$ is a irreducible component of either $\Do^\ast$ or $\Dab^\ast$. Note that in any case $D^{\ast,c} \cong S(p^2)$, and by Yamazaki's results
\[ \begin{array}{lcl}
L^2 D^{\ast,c} &=& 0 \\
L (D^{\ast,c})^2 &=& - \frac{p^6(p^2-1)}{12}. 
\end{array} \]
Thus the trace formula reads
\[ \trace(\gamma^\ast | S_k(\Gamma_2(p^2))) = - k \cdot \frac{\xi^{-1}}{(1- \xi^{-1})^2}  \: L (X^\gamma)^2 + \ldots \]
if $\gamma$ is a generator of the stabilizer of $D^{\ast,c}$. 
To get the total contribution of the boundary components we have to add up the contributions of all boundary components involved, and we have to consider the sum over all nontrivial $p^2$-th roots of unity. The following lemma is not difficult to prove.

\begin{lemma}\label{7}
Let $\xi$ be a $n$-th root of unity and primitive. Then
\[ \sum_{j=1}^{n-1} \frac{1}{1-\xi^j} = \frac{n-1}{2} \]
as well as
\[ \sum_{j=1}^{n-1} \frac{\xi^i}{(1-\xi^i)^2} = - \frac{n^2-1}{12} .\]
\end{lemma}

Applying the same considerations as in case 1.a one finds that there are
\[ \frac{p^5(p^2-1)}{2} + p^3 \, \frac{p(p-1)}{2} = \frac{1}{2} \, p^4(p^3-1)\]
irreducible components in all of $\Do^\ast$ and $\Dab^\ast$. Hence the contribution from the boundary is
\[ \begin{array}{lcl}
t_3 &=& -k \cdot (- \frac{p^4-1}{12}) \frac{p^4(p^3-1)}{2} ( - \frac{p^6(p^2-1)}{12}) + \ldots\\[10pt]
&=& - k \cdot \frac{p^{10}(p^2-1)^2}{12 \cdot 24} (p^2+1)(p^3-1) + \ldots\\[10pt]
&=& - k \cdot \frac{p^{10}(p^2-1)^2}{2^8 3^3 5} \, 120 (p^5 + p^3 - p^2 - 1) + \ldots
\end{array} \]

2) Finally we have to consider the fixed varieties of codimension 2. In these cases the trace formula becomes
\[ \begin{array}{l}
\trace(\gamma^\ast | S_k(\Gamma_2(p^2))) =\qquad \\[2mm]
\quad =  \frac{1}{1-e^{{-i \Theta_1}}} \: \frac{1}{1-e^{{-i\Theta_2}}} \: ( 1+ c_1(\calSk | X^\gamma) )  \cdot ( 1 - \frac{1}{2} K_{{X^\gamma}}) \\[2mm]
\qquad \cdot ( 1- \frac{e^{{-i\Theta_1}}}{1-e^{{-i\Theta_1}}} c_1(N^\gamma(\Theta_1))) \cdot  
 ( 1- \frac{e^{{-i\Theta_2}}}{1-e^{{-i\Theta_2}}} c_1(N^\gamma(\Theta_2))) 
\\[2mm]
\quad = k^1 \cdot \frac{1}{1-e^{{-i\Theta_1}}} \frac{1}{1-e^{{-i\Theta_2}}} \: \mbox{deg} \, L | X^\gamma \: + \: k^0 \cdot ( \ldots) .
\end{array} \]
Here we assumed $\chi(\gamma) = 1$, see equation $(\dagger)$ above. In case 2.c below we will have to change the sign of this expression, since then $\chi(\gamma) = -1$. 
   
\underline{Case 2.a:} Let $\gamma = \zeta_1$ and $X^\gamma = C_1^\ast$ with $e^{{i \Theta_1}} = -1$ and $e^{{i \Theta_2}} = -i$.
There are precisely $\frac{p^{10}(p^2-1)}{4}$ components of $C_1^\ast$. Each one is isomorphic to $X(p^2)$. Summing up all elements of the stabilizer which are not equal to $\Eins$ or $\nu_1$ and using lemma \ref{7}  we get as total contribution from $C_1^\ast$
\[ t_4 = k \cdot \frac{1}{2} \frac{p^4(p^2-1)}{24} \frac{p^{10}(p^2-1)}{4}  + \ldots \]

\underline{Case 2.b:} Let  $\gamma = \zeta_2$ and $X^\gamma = C_2^\ast$ with $e^{{i \Theta_1}} = -\rho^2$ and $e^{{i \Theta_2}} = - \rho$ where $\rho = e^{2\pi i / 3}$.
An entirely analogous calculation shows that the total contribution from the $\frac{p^{10}(p^2-1)}{6}$  components of $C_2^\ast$ is
\[ t_5 = k \cdot \frac{2}{3}  \frac{p^4(p^2-1)}{24} \frac{p^{10}(p^2-1)}{6}   + \ldots \]
As combined contribution we get
\[ t_4 + t_5 = k \cdot  \frac{(p^2-1)^2}{2^8 3^3 5} \: 340 p^{14}  + \ldots \]

\underline{Case 2.c:} The intersection curves $B_i^\ast$ for $i=1,2$. Let $B_i^{\ast,c}$ be a component of $B_i^\ast$ with stabilizer of order $2p^2$, generated by $\nu_i$ and the element of order $p^2$ that generates the stabilizer of the corresponding boundary component. This gives $e^{{i \Theta_1}} 
= -1$ and $e^{{i\Theta_2}} = e^{{2\pi i / p^2}}$.
Note however that in this case we have 
\[ \mbox{ch}\,(\calSk | B_i^{\ast,c})(\gamma) = - \mbox{ch}\,(\calSk | B_i^{\ast,c}) \]
for any $\gamma$ in the stabilizer of $B_i^{\ast,c}$ which is not in the stabilizer of the boundary component. 
This is because the action of such a $\gamma$ on $\oka_{{\alppast}}(-D)|B_i^{\ast,c}$ which is isomorphic to the conormal bundle of $B_i^{\ast,c}$ in $H_i^\ast$  is not trivial, but a reflection on the fibres. 
The degree of the map $\Do^{\ast,c} \rightarrow \Do$ is $2 p^6$ with ramification of order $2$ above $\Do \cap H_i$, $i=1,2$ as shown in \cite{Z}. Since the degree of $X(p^2) \rightarrow X(p)$ equals $p^3$ there are $p^3$ components of $B_1^\ast$ in each component of $\Do^\ast$. The degree of $\Dab^{\ast,c} \rightarrow \Dab$ for a nonstandard component $\Dab$ equals $p^8(p^2-1)$ with ramification of order $2$ as before, and the degree of $X(p^2) \rightarrow X(1)$ is $\frac{p^4(p^2-1)}{2}$. Hence we get $p^4$ components of $B_1^\ast$ in each $\Dab^{\ast,c}$. Thus there are precisely
\[ p^4 \cdot p^3 \cdot \frac{p(p-1)}{2} + p^3 \cdot \frac{p^5(p^2-1)}{2} = \frac{p^8}{2} (p^2 + p - 2) \]
irreducible components of $B_1^\ast$. Summing up over all elements of the stabilizer, which were not already contained in the stabilizer of either $H_1^\ast$ or the boundary component, we get 
\[ \begin{array}{lcl}
t_6 &=& - k^1 \cdot \frac{p^4(p^2-1)}{24} \frac{p^8(p^2+p-2)}{2} \frac{1}{2} \frac{p^2-1}{2} + k^0 \cdot (\ldots) \\[2mm]
&=& - k^1 \cdot \frac{(p^2-1)^2}{2^8 3^3 5} \cdot  180 ( p^{14} + p^{13} - 2 p^{12} ) + \ldots 
\end{array}  \]
The contribution of $B_2^\ast$ is three times the contribution of $B_1^\ast$ since the $6:1$ covering map $H_2^{\ast,c} \rightarrow H_1^{\ast,c}$ is ramified of order two along the components of $B_2^\ast$. Hence
\[ t_7 = 3 \cdot t_6 .\]

Summing up all contributions we get from equation $(\ast)$
\[ \begin{array}{l}
d \cdot \dim S_k(\gap) - \dim S_k(\Gamma_2(p^2)) = t_1 + \ldots + t_7 =\\[10pt]
= \frac{(p^2-1)^2}{2^8 3^3 5} [ \: k^2 \: \cdot 210 p^{14} \\[10pt]
\qquad \qquad + k^1 \cdot (-120 p^{15} -1010 p^{14} - 840 p^{13}+120 p^{12} + 120 p^{10} )\\[10pt]
\qquad \qquad + k^0 \cdot ( \ldots) ] .
\end{array} \]
Now compare this with the dimension formulas obtained via Riemann--Roch:
\[ \begin{array}{l}
d \cdot \dim S_k(\gap) - \dim S_k(\Gamma_2(p^2)) = \\[10pt]
= \frac{(p^2-1)^2}{2^8 3^3 5} [ \: k^2 \cdot 210 p^{14} \\[10pt]
\qquad \qquad + k^1 \cdot (-4 p^{16} -120 p^{15} - 1010 p^{14} + 120  p^{13}+120 p^{12} + 120 p^{10} )\\[10pt]
\qquad \qquad + k^0 \cdot ( \ldots) ]\\[10pt]
\qquad + k^1 \cdot  \frac{1}{12} \, d \cdot c_2(\alptilde).L \: . 
\end{array} \]
Note that we don't need to know the numbers $c_1.c_2$ and $c_2(\alptilde).\Do$ since they only contribute to $k^0$. 
This immediately yields the following result.

\begin{prop}
The intersection number $c_2(\alptilde).L$ equals
\[ c_2(\alptilde).L = \frac{p^2-1}{720} (p^3 + 121 p + 60 ). \]
\end{prop}

The constant term of this formula differs from the one conjectured in \cite{Z}. This is because a mistake was made there in the computation of the contribution from $c_2(\alptilde).\Do$. \\

However, this does not affect the conjecture on the second Chern number. In \cite{Z} we conjectured that $c_1.c_2$ would be a polynomial in $p$ of certain type, whose value for small $p$ is known. Hence we could predict a formula, which we will prove now.

\begin{thm}
The second Chern number of $\alptilde$ equals
\[ c_1(\alptilde).c_2(\alptilde) = - \frac{p^2-1}{240} (p-13)(p^2-17p+90) . \]
\end{thm}

\proof
Using theorem \ref{1} the second Chern number can be calculated by computing $c_2(\alptilde).L$ and the product of $c_2(\alptilde)$ with $H_1$, $H_2$, $\Do$, $\Dab$ and the exceptional divisors.
Note that apart from $\Do$ all of these are smooth hypersurfaces in $\alptilde$. Therefore we have
\[ c_2(\alptilde).H_1 = c_2(H_1) - c_1(K_{{H_1}}).c_1(\calN_{{H_1/\alptilde}}) \] 
where $c_2(H_1)$ is equal to the Euler number of $H_1$, 
and similarly for the other intersection numbers. Using the description of the hypersurfaces in \cite{HKW} one can verify that these numbers are precisely those given  in table 4. See \cite{Z} for details. 

 \pagebreak

{\em Table 4: }
{Euler numbers and intersection numbers for certain smooth hypersurfaces $A \subset \alptilde$}
\[ \begin{array}{l||r|r}
A&\mbox{e}\,(A)&c_1(K_A).c_1(\calN_{{A/\alptilde}}) \\
\hline
H_1& - \frac{p-6}{6}(p^2-1) & \frac{p+6}{12}(p^2-1)\\[1pt]
H_2& - \frac{p-21}{6}(p^2-1) & \frac{p+3}{12}(p^2-1)\\[1pt]
\Dab&9&1\\[1pt]
E^{(1)}& - \frac{p-6}{6}(p^2-1) & p^2-1\\[1pt]
E^{(2)}& - \frac{p-6}{6}(p^2-1) & \frac{3}{2}(p^2-1)\\[1pt]
E^{(3)}_{a,b}&3&6\\[1pt] 
E^{(4)}_{a,b}& 4&6\\[1pt]
E^{(5)}_{a,b}& 3&6
\end{array} \]

 The number $c_2(\alptilde).L$ was determined above. However, the constant term of the formular for $c_2(\alptilde).\Do$ stated in \cite[lemma 6.4.8]{Z} is not correct.
The difficulty in computing this particular number arises from the fact that the central boundary component $\Do$ is not a smooth surface. \\

But we know that $c_1.c_2$ is a birational invariant of the moduli space. So we can as well compute the second Chern number by looking at a smooth blow up $\mu : \alpr \rightarrow \alptilde$ which resolves all of the singularities of $\Do$. Such a resolution was constructed explicitely in \cite[section 6.4]{Z}, where it was denoted by $X''$. First, $\alptilde$ has to be blown up in all {\em inner deepest points}, leading to an exceptional divisor $G''$, which consists of $\frac{(p^2-1)(p-5)}{12}$ pairwise disjoint components. Then we blow up along $\frac{(p^2-1)(p-3)}{8}$ pairwise disjoint {\em inner rational curves}, giving an exceptional divisor $F''$. \\

The strict transform of a divisor $A$ on $\alptilde$ under $\mu$ shall be denoted by $A''$. We get
\[ \begin{array}{lll}
K_{{\alpr}} &=& \mu^\ast K_{\alptilde} + F'' + 2 G'' \\
&=& 3 L'' - \Do'' - \sum \Dab'' - \frac{1}{2} H_1'' - \frac{1}{2} H_2'' 
 - E'' - F'' - G'' 
\end{array} \]
because
\[ \mu^\ast \Do = \Do'' + 2 F'' + 3 G'' .\]
Here $E''$ stands for the contribution from the exceptional divisors on $\alptilde$, see theorem \ref{1}. 
Since for a smooth divisor $A''$ on $\alpr$ one has
\[ \begin{array}{lll}
c_2(\alpr).A'' &=& c_2(A'') - K_{{A''}} . A''|A''\\
&=& c_2(A'') - K_{{A''}}.(K_{{A''}}- K_{{\alpr}}|A'') 
\end{array} \]
we have the equality
\[ c_2(\alpr).\mu^\ast A = c_2(\alptilde).A \]
if a hypersurface $A$ on $\alptilde$  does not contain any of the points,  which get blown up under $\mu$. 
Since $L$ restricted to the boundary of $\alptilde$ is trivial and $\mu$ blows up points of the boundary only, we also have
\[ c_2(\alpr).L'' = c_2(\alptilde).L .\]
By \cite{Z} or \cite{HKW} the only hypersurfaces from the above affected by $\mu$ are $H_2$ and $\Do$. Hence to compute $c_1(\alpr).c_2(\alpr)$ we have to compute $c_2(\alpr).H_2''$, $c_2(\alpr).\Do''$, $c_2(\alpr).F''$ and $c_2(\alpr).G''$. All other numbers needed are equal to the numbers given for $\alptilde$ in table 4. In fact, in \cite[Hilfssatz 6.4.5]{Z} we already computed
\[ c_2(\alpr).\Do'' = \frac{p^2-1}{24} ( 3p^2 - 10 p + 3) + 6 \, \frac{(p^2-1)(p-5)}{12} .\]
$(i)$ Let us look at $H_2''$ first. This Humbert surface contains precisely $\frac{p^2-1}{2}$ points on inner rational curves, which are deepest points, but no inner deepest points. It is smooth in each of this points but of multiplicity $2$ along the inner rational curve. Hence if $F''_{{H_2}}$ denotes the reduced exceptional divisor of the blow up in $H_2''$ we get $F''|H_2'' = 2 F_{{H_2}}''$ and $(F''_{{H_2}})^2 = -1$. We have
\[ \begin{array}{llll}
&K_{{H_2''}} &=& \mu^\ast K_{{H_2}} + F''_{{H_2}} \\
\mbox{and} \:\: &H_2''|H_2'' &=& \mu^\ast H_2|H_2 - 2 F''_{{H_2}} .
\end{array} \] 
Now we compute
\[ \begin{array}{lll}
c_2(\alpr).H_2'' &=& c_2(H_2'') - K_{{H_2''}} . H_2''|H_2''\\
&=& c_2(H_2) + \frac{p^2-1}{2} - (\mu^\ast K_{{H_2}} + F''_{{H_2}}).(\mu^\ast H_2|H_2 - 2 F''_{{H_2}}) \\
&=& c_2(H_2) + \frac{p^2-1}{2} - K_{{H_2}}.H_2|H_2 + 2 (F''_{{H_2}})^2\\
&=& \frac{p^2-1}{12} (-3p + 33) 
\end{array} \] 
using the fomulas above and the numbers from table 4. \\

$(ii)$ Next we consider $G''$. By construction a component $G^c$ of $G''$ is isomorphic to $\PP^2$ blown up in three points. Let Pic$(G'')$ be generated by the class of a general line $h$ and the exceptional divisors $\varepsilon_1$, $\varepsilon_2$ and $\varepsilon_3$. We find
\[ K_{{G^c}} = -3 h + \varepsilon_1 + \varepsilon_2 + \varepsilon_3 \]
and
\[ G^c|G^c = K_{{G^c}} - K_{{\alpr}}|G^c = -h .\]
Hence we obtain
\[ \begin{array}{lll}
c_2(\alpr).G'' &=& \frac{(p^2-1)(p-5)}{12} \: ( c_2(G^c) - K_{{G^c}}.G^c|G^c ) \\[6pt]
&=& 3 \cdot \frac{(p^2-1)(p-5)}{12} .
\end{array} \]
$(iii)$ Finally we have to deal with $F''$. There are two types of components $F^c$ of $F''$: Those which are top components and those which are not. See \cite[section 6.4]{Z} for details. \\

The $\frac{p^2-1}{2}$ top components are isomorphic to $\PP(\oka_{{\PP^1}} \oplus \oka_{{\PP^1}}(1) )$ with
\[ \begin{array}{lll}
K_{F^c} &=& -2 b - 3 f \\
F^c|F^c &=& -b -3f 
\end{array} \]
where $b$ and $f$ are as before.
Therefore we get 
\[ c_2(\alpr).F^c = c_2(F^c) - K_{{F^c}}.F^c|F^c = -3 .\]
The remaining $\frac{p^2-1}{2}(\frac{p-3}{4} -1)$ components $F^c$ are isomorphic to $\PP^1 \times \PP^1 $ and we have
\[ \begin{array}{lll}
K_{F^c} &=& -2 b - 2 f \\
F^c|F^c &=& -b -4f .
\end{array} \]
Here $b$ and $f$ denote the usual generators of Pic$(\PP^1 \times \PP^1)$. 
Hence we compute
\[ c_2(\alpr).F^c = -6 \]
and taking everything together we obtain for $F''$
\[ c_2(\alpr).F'' = \frac{p^2-1}{8}( -6p + 30) .\]
Now we have computed all of the contributions to $c_1(\alpr).c_2(\alpr)$. Adding up and using the results from table 4 we finally obtain the formula as claimed. 
\ebew

\begin{cor}
The arithmetic genus of the moduli space $\alptilde$ equals
\[ p_a(\alptilde) = 1 - \frac{c_1.c_2}{24} = 1 + \frac{p^2-1}{5760} \: (p-13)(p^2-17p+90) . \]
\end{cor}

\begin{rk}\em
If $p=5,7$ or $11$ then $p_a(\alptilde) = 0$, because in each of these cases $\alptilde$ is unirational, see \cite{MS,GP}. Since we have  $p_a(\alptilde) \neq 0$ for $p \ge 13$ none of these moduli spaces can be unirational, in accordance with Gritsenko's result \cite{Gr}.
\end{rk}

In retrospect we can now state the correct formula for $c_2(\alptilde).\Do$ replacing lemma 6.4.8 of \cite{Z}.

\begin{lemma}
The intersection number $c_2(\alptilde).\Do$ equals
\[ c_2(\alptilde).\Do = \frac{p^2-1}{24} ( 3p^2 - 10p - 3) . \]
\end{lemma}

An immediate consequence of our computations using the theorem of \linebreak Riemann--Roch is a polynomial dimension formula for cusp forms with respect to $\gap$. There are several results on dimension formulas which are more general, but usually not as explicit as ours.  See for example \cite{Has}, which also contains further references. 

\begin{thm} 
Let $k \neq 0$ and $k \equiv 0$ mod $12$. Then the space of cusp forms of weight $k$ with respect to $\gap$ has dimension
\[ \begin{array}{crcl}
\dim S_k(\gap)\hspace{5mm} =
& \frac{p^2-1}{2^8 3^3 5} \cdot [ \: \: \:  & k^3  \cdot &(2 p^{3} + 2 p)\\
&+&k^2  \cdot & (-9 p^{3} + 201  p)\\
&+&k^1 \cdot & (13 p^{3} - 120  p^{2} -997 p - 840  )\\
&+&k^0 \cdot & (-6)  (p-35) (p+2) (p+3) \: \: ] .
\end{array} \]
\end{thm}

\noindent
J\"org Zintl\\
Fachbereich Mathematik\\
Universit\"at Kaiserslautern\\
67653 Kaiserslautern\\
Germany\\
{\em zintl@mathematik.uni-kl.de}


\begin{thebibliography}{HKW 2}


\bibitem[AMRT]{AMRT}{Ash A., Mumford D., Rapoport M., Tai Y., {\em Smooth Compactification of Locally Symmetric Varieties}, Math Sci Press, Brookline (1975)}

\bibitem[AS]{AS}{Atiyah M., Singer I., The index of elliptic operators: III, {\em Ann. of Math.} 87 (1968) 546--604  }
	
\bibitem[G]{Gr}{Gritsenko V., Irrationality of the moduli spaces of polarized abelian surfaces, in: Barth et al., {\em Abelian Varieties - Proceedins of the Egloffstein conference} (1993), de Gruyter, Berlin (1995) 63--81}   

\bibitem[GH]{GH}{Gritsenko V., Hulek K., Appendix to the paper ``Irrationality of the moduli spaces of polarized abelian surfaces'', in: Barth et al., {\em Abelian Varieties - Proceedins of the Egloffstein conference} (1993), de Gruyter, Berlin (1995) 83--84}

\bibitem[GP]{GP}{Gross M., Popescu S., Calabi Yau threefolds and moduli of abelian surfaces, announced at the Santa Cruz conference (1995)}

\bibitem[Ham]{Ham}{Hammond W.F., On the graded ring of Siegel modular forms of genus two, {\em Am. J. Math.} 87 (1965) 502--506}

\bibitem[Har]{Harder}{Harder G., A Gauss--Bonet formula for discrete arithmetically defined groups, {\em Ann. Sci. Ecole Norm. Sup.} (4) 4 (1971) 409--455}

\bibitem[Has]{Has}{Hashimoto K., The dimension of the spaces of cusp forms on Siegel upper half plane of degree two (I), {\em J. Fac. Sci. Univ. Tokyo} 30 (1983) 403--488} 

\bibitem[Hi]{Hi}{Hirzebruch F., Elliptische Differentialoperatoren auf Mannigfaltigkeiten, in: {\em Gesammelte Abhandlungen / Collected Papers}, Vol. II, 37,  Springer, Berlin  (1987)}

\bibitem[HM]{HM}{Horrocks G., Mumford D., A rank 2 vector bundle on $\PP^4$ with 15,000 symmetries, {\em Topology} 12 (1973) 63--81}

\bibitem[HKW 1]{HKWsing}{Hulek K., Kahn C., Weintraub S.H., Singularities of the moduli spaces of certain abelian surfaces, {\em Compositio Math.} 79 (1991) 231--253}

\bibitem[HKW 2]{HKW}{Hulek K., Kahn C., Weintraub S.H., {\em Moduli spaces of Abelian Surfaces: Compactification, Degenerations and Theta Functions}, de Gruyter, Berlin (1993)}

\bibitem[HS 1]{HS}{Hulek K., Sankaran G.K., The Kodaira dimension of certain moduli spaces of abelian surfaces, {\em Compositio Math.} 90 (1994) 1--36}

\bibitem[HS 2]{HS2}{Hulek K., Sankaran G.K., The fundamental group of some Siegel modular threefolds, in: Barth et al., {\em Abelian Varieties - Proceedings of the Egloffstein conference} (1993), de Gruyter, Berlin (1995) 141--150}

\bibitem[I]{I}{Igusa J., On Siegel modular forms of genus two (II), {\em Am. J. Math.} 86 (1964) 392--412}

\bibitem[K]{K}{Kawamata Y., A generalization of Kodaira--Ramanujam's vanishing theorem, {\em Math. Ann.} 261 (1982) 43--46}

\bibitem[MS]{MS}{Manolache N., Schreyer F.-O., Moduli of (1,7)-polarized abelian surfaces via syzygies, {\em preprint} math.AG/9812121}

\bibitem[Sa]{Sa}{Satake I., On the compactification of the Siegel space, {\em J. Indian Math. Soc.} 20 (1956) 259--281}

\bibitem[Y]{Y}{Yamazaki T., On Siegel modular forms of degree two, {\em Am. J. Math.} 98 (1976) 39--53}

\bibitem[Z]{Z}{Zintl J., {\em Invarianten von kompaktifizierten Modulr\"aumen polarisierter abelscher Fl\"achen}, Thesis, Hannover (1996), available at: http://www.mathematik.uni-kl.de/$\sim$wwwagag/E/Zintl/}

\end{thebibliography}
\end{document}